\documentstyle[12pt]{article}

\newcommand{\dfrac}{\displaystyle\frac}
\newcommand{\gru}[2]{\mathop{{\bf #1 }(#2)}\nolimits}
\newcommand{\ali}[2]{\mathop{{\sl #1 }\,(#2)}\nolimits}

\newcommand{\ma}[4]{\mbox{$\pmatrix{#1&#2\cr#3&#4}$}}

\newtheorem{Prop}{Proposition}
\newtheorem{Theo}{Theorem}
 
\newtheorem{Definition}{Definition}
\newtheorem{ex}{Example}
\newenvironment{Example}{\begin{ex}\rm}{\end{ex}}

\normalbaselines
\begin{document}
\begin{center}
\vspace*{1.0cm}

{\LARGE{\bf  Class Operators as Intertwining Maps into the Group Algebra}}

\vskip 1.5cm

{\large {\bf Aleksander Strasburger}}

\vskip 0.5 cm

Department of Mathematical Methods of Physics,\\
University  of Warsaw\\
Ho\.{z}a 74,\\
00-682 Warszawa, Poland \\
e-mail: \verb+alekstra@fuw.edu.pl+
\end{center}

\vspace{1 cm}

\begin{abstract}
With the aim of completing the previous study by  A. Or{\l}owski and the author concerning 
intertwining maps between induced representations and conjugation representation, termed here weighted class operators, we compute the latter explicitely for the conjugation representation arising from the regular representation in the group algebra of a compact group. To that efect a theorem of Wigner--Eckart type for weighted class operators obtained from matrix coefficients of irreducible representations of a compact group is proved. Also the previous construction of weighted class operators is reviewed and extended to the case of locally compact groups rather then just compact ones.

\end{abstract}

\vspace{1 cm}

\section{Introduction}

\footnotetext[1]{1991 {\em Mathematics Subject Classification}. Primary 22A25, 22D15, 22D30, Secondary 20C35} 

In recent papers \cite{OS,SO} A. Or{\l}owski and the author have
investigated, in the context of the theory of representations of compact
groups, a certain construction extending naturally that of the class
operator. Let us remind that the class operators, as they are defined in
the context of the theory of finite groups, are elements of the group
algebra which are sums of group elements belonging to a particular
conjugacy class. They span an abelian subalgebra of the group algebra and
the knowledge of this subalgebra (including its multiplicative structure)
is tantamount to the knowledge of the representation theory of the group
(see \cite{Ch} and also the begining of the Section 3 of the present
paper, where we elaborate slightly on this point). The extention alluded
to above consists in replacing class sums by operators obtained
by integrating the class functions against conjugation operators of the
representation --- the class operator corresponds to the class function
being identically 1. Here the term class function describes a function
whose domain is a fixed conjugacy class. In view of an interpretation of
such functions as (non-necessarily positive) weights assigned to points of
the conjugacy class, the operators obtained this way were called in
\cite{SO} weighted class operators. 

In the mentioned papers some properties of this construction, focusing on
properties of the map assigning weighted class operators to class
functions, were investigated for compact groups. However since the concept
is expected to play a significant role for the non-compact case as well,
we describe here the relevant construction for locally compact rather then
just for compact groups. 
     
\section{The construction of weighted class operators} 

Given any locally compact topological space $X$ we shall denote by ${\cal
K}(X)$ the space of continuous compactly supported functions on $X$,
endowed with the topology of uniform convergence on compact subsets. 

If $G$ is a unimodular locally compact group and $V$ a topological vector
space then by a representation of $G$ on $V$ we shall mean a homomorphism
$T\colon G\to GL(V)$ into the group $GL(V)$ of continuous invertible
linear maps on $V$, which is continuous with respect to the strong
operator topology. In such a case we shall simply say $(T,\,V)$ is a
representation of $G$. In the most cases we shall deal with unitary
representations (not necessarily irreducible) of $G$ on a Hilbert space
$V$ by what we shall mean a homomorphism $T\colon G\to U(V)$, where
$U(V)\subset GL(V)$ is the group of unitary automorphisms of $V$. 

If $L(V)$ denotes the space of continuous linear operators on $V$ with the
strong topology, then the conjugation $L(V)\ni A\mapsto T(g)AT(g^{-1})\in
L(V)$ defines a continuous representation of $G$ on $L(V)$ which will be
called the conjugation representation defined by $(T,\,V)$. 

If $f\in {\cal K}(X)$ and $g_0\in G$ is an arbitrary, but fixed, element
of the group, then we define
\begin{equation}
 T(f;\,g_0)= \int_{G} f(x)T(x)T(g_0)T(x^{-1})\,dx \label{int}
\end{equation} where the integration is performed with respect to the
Haar measure $dx$ on $G$. It can easily be shown by the
standard methods (cf. e.g. \cite{AR,AW}) that the integral (\ref{int})
converges in the sense of the strong topology in $L(V)$. 
Thus for each fixed $g_0\in G$ we have a linear mapping
\begin{equation}\label{def}
{\cal K}(G)\ni f\mapsto   T(f;\,g_0)\in L(V). 
\end{equation}

By $\lambda$, resp.~$\rho$, we shall denote the left, resp.~right, regular
representation of $G$ in ${\cal K}(G)$, where the action of an element
$g\in G$ is defined as the mapping $f\mapsto\lambda(g)f$, resp.  $f\mapsto
\rho(g)f$, where 
\begin{equation}\label{regular}
\lambda(g)f(x)=f(g^{-1}x), \quad\mbox{resp.}\quad \rho(g)f(x)=f(xg),\quad
x\in G. 
\end{equation} 
A straightforward computation gives the following relations 
\begin{eqnarray} 
T(g)T(f;g_{0})T(g)^{-1}& = & T(\lambda(g)f;g_{0}), \label{fun}\\ 
T(\rho (h)f;g_{0})& =& T(f;g_{0}),\qquad
\mbox{\rm for each $h\in Z_{0}$}, \label{quo} 
\end{eqnarray} 
where $Z_{0}=Z(g_{0})\subset G$ denotes the centralizer of $g_{0}$ in $G$
defined by $h\in Z_{0} \iff hg_{0}=g_{0}h$. Note $Z_{0}$ is a closed
subgroup of $G$. 

Recall now that the conjugacy classes in $G$ can be regarded as
$G$-homogeneous spaces in the following way. For any element $g_0 \in G$
let $C_0= C(g_0)=\{xg_0x^{-1}\bigm|x\in G\}$ be its conjugacy class in
$G$. The map $G\ni x \mapsto xg_0x^{-1} \in C_ 0$ is surjective and
constant on the left $Z_0$ cosets in $G$, hence it induces a bijection
$G/Z_0\ni xZ_0\mapsto xg_0x^{-1}\in C_0$ of $G/Z_0$ with the conjugacy
class $C_0$, such that the left action of $G$ on $G/Z_0$ corresponds to
the action by conjugation on $C_0$.  We shall use this bijection to
transfer topology and in the Lie case also the manifold structure from
$G/Z_0$ to the conjugacy class $C_0$. The reason for doing so is that the
subspace topology inherited by $C_0$ from $G$ in general is rather
complicate d, in particular $C_0$ need not be closed in $G$, while the
topology, and if aplicable also the manifold structure of the homogeneous
space $G/Z_0$ are much simpler. However, if G is compact, then all its
conjugacy classes are compact subsets of $G$ and the topology transferred
from $G/Z_0$ by means of the above map is the same as their subspace
topology. 
  
The invariance condition (\ref{quo}) allows us now to use the following
well known technique (cf. e.g. \cite{AW}) of transferring integrals from
the group $G$ to the homogeneous space $G/Z_0$, thus passing from the map
(\ref{def}) to the map of the space of functions on the coset space
$G/Z_0$.  The transfer can be described as follows.  Let $dh$ denote the
(left) Haar measure on $Z_0$ and assume for simplicity that there exists
an invariant measure $d\mu$ on $G/Z_0 $ such that \begin{equation}
\int_Gf(x)\,dx=\int_{G/Z_0 }d\mu(\dot{x})\int_{Z_0 }f(xh)\,dh,
\label{InvInt} \end{equation} for each $f\in {\cal K}(G)$, where the
symbol $\dot{x}$ is used to denote the coset $xZ_0 \in G/Z_0 $
corresponding to $x\in G$.  It is known (cf. e.g. \cite{AW}) that the
mapping obtained by averaging functions on $G$ over $Z_0$-cosets, to wit 
$f\mapsto \widetilde{f}$, where
\begin{equation}
\widetilde{f}(\dot{x})=\int_{Z_0}f(xh)\,dh, \label{transfer}
\end{equation} 
is a linear surjection of ${\cal K}(G)$ onto ${\cal
K}(G/Z_0)$, carrying the left regular representation of $G$ in ${\cal
K}(G)$ onto the natural representation of $G$ by left translations in
${\cal K}(G/Z_0)$.  We shall therefore denote the latter also by $\lambda$
and note for the future use that by the natural extention by continuity
(cf. below) it gives rise to the representation of $G$ induced by the
trivial (one dimensional identity) representation of $Z_0$.  Now observing
that $T(xg_0x^{-1})$ depends only on the coset $\dot{x}=xZ_0 $ of $x$, so
we can write $T(xg_0x^{-1})=T(\dot{x})$, we use (\ref{InvInt}) to rewrite
the equation (\ref{int}) in the form 
\begin{eqnarray*}
 T(f;\,g_0)& =& \int_{G} f(x)T(x)T(g_0)T(x^{-1})\,dx =
\int_{G/Z_0}d\mu(\dot{x})T(\dot{x})\int_{Z_0 }f(xh)\,dh\\
 & =&\int_{G/Z_0 }\widetilde{f}(\dot{x})T(\dot{x})\,d\mu(\dot{x}). 
\end{eqnarray*}
This shows that the mapping (\ref{def}) factorises through the surjection
${\cal K}(G)\ni f\mapsto \widetilde{f}\in {\cal K}(G/Z_0)$, giving rise to
a map ${\cal K}(G/Z_0 )\to L(V)$. We shall collect properties of this
construction in the following form.  
\begin{Prop}{\rm(cf. \cite{OS,SO})}
\label{Propos1} 
Let for $f\in {\cal K}(G)$ the operator $T(f;\,g_0)$ be
defined by {\rm(\ref{int})} and for $\phi\in {\cal K}(G/Z_0 )$ define the
operator $\widetilde{T}(\phi;\,g_0)\in L(V)$ by setting
\begin{equation}
\widetilde{T}(\phi;\,g_0)  = 
\int_{G/Z_0 }\phi(\dot{x})T(\dot{x})\,d\mu(\dot{x}).  
 \label{int2} \end{equation} Then the mapping 
\begin{equation}
\label{def2} 
\widetilde{T}\colon {\cal K}(G/Z_0 )\ni \phi \mapsto
\widetilde{T}(\phi;\,g_0) \in L(V) \end{equation} satisfies the condition
of covariance with respect to the action of $G$,
 \begin{equation}
T(g)\widetilde{T}(\phi;\,g_{0})T(g)^{-1}=
\widetilde{T}(\lambda(g)\phi;\,g_{0}), \quad \phi\in {\cal K}(G/Z_0),
\label{imprim} 
\end{equation} 
and for every $f\in {\cal K}(G)$ such that
$\widetilde{f}= \phi$ we have
\begin{equation} T(f;\,g_0)=
\widetilde{T}(\phi;\,g_0) .  
\end{equation} 
In addition, if the representation $(T,V)$ is unitary, then the map
{\rm(\ref{def2})} extends by continuity to the mapping of the space
$L^2(G/Z_0) $ of all (equivalence classes of) square integrable functions
on $G/Z_0$ to the space of Hilbert--Schmidt operators on $V$. 
\end{Prop}

The maps (\ref{def}), (\ref{def2}) are the generalizations of the class
operator we have put forward (in the case of the compact group) in
\cite{OS,SO}. We shall refer to them indifferently as a weighted class
operator maps for the representation $(T,V)$ based on the conjugacy class
$C_0=G/Z_0$.  Note, however that the class operator, which should
correspond to the integral $\int_{G/Z_0 }T(\dot{x})\,d\mu(\dot{x})$ might
be not defined itself, unless the group is compact or the representation
is rather special, since in general this integral does not converge. 

It is clear that the covariance condition (\ref{imprim}) is implied by the
equation (\ref{fun}).  One also says the map (\ref{def2}) intertwines the
 action of $G$ by $\lambda$ on ${\cal K}(G/Z_0)$ with the conjugation
representation on $ L(V)$.  The significance of the last part of the
statement of Proposition 1 can be grasped better, if one notes the
representation $(\lambda,\,L^2(G/Z_0) )$ of $G$ is nothing else but the
representation of $G$ induced by the trivial representation of $Z_0$. 
Thus we see that the construction gives an intertwining operator between the
induced representation and the conjugation representation. We shall say
more to this point in the  last section. 

We point out the double sided character of the Proposition \ref{Propos1}.
On one side the equation (\ref{int2}) is closer to the original meaning of
the class operator, while on the other hand the actual computations are
often easier to handle on the level of the group, by the use of the
formula (\ref{int}), as we shall see below. 

\begin{Example}

We  illustrate the use  of the above technique of integration on a classical example of the group $\gru{SU}{2}$, in order to show that the equation (\ref{int2}) for the function $\phi$ being identically  $1$, leads  in that case to the integral 
\begin{equation}
\frac{1}{4\pi}\int^{2\pi}_0\int^\pi_0
\exp[i\psi(J_x\sin\theta\cos\varphi+
J_y\sin\theta\sin\varphi+
J_z\cos\theta)]\sin\theta\,d\theta\,d\varphi, \label{FR}
\end{equation}
whose evaluation was the main issue in the papers  \cite{FR,B,MZ} (also cf.  \cite{OS}). 
The  $J_x,\,J_y,\,J_z$ are of course the infinitesimal generators of a
representation of the group $\gru{SU}{2}$ (or the rotation group $\gru{SO}{3}$).  In the present formalism it is clearly enough to consider the case of the identity (spin $1/2$) matrix representation of  $\gru{SU}{2}$, 
\[
\gru{SU}{2} =
\Biggl\{\ma{\hphantom{-}a}{b}{-\overline{b}}{\overline{a}}\biggm|a,\,b\in{\bf C},\ 
\vert a\vert^2 + \vert b\vert^2 = 1\Biggr\},
\]
in which case the $J $'s will be replaced by the Pauli matrices 
$\sigma_\alpha$.  
The Lie algebra $\ali{su}{2}$ of $\gru{SU}{2}$
is therefore  the space of antihermitean traceless $2\times2$ matrices, which can be identified with  ${\bf R}^3$ by means of the basis   $\{i\sigma_\alpha\}_{\alpha=1}^3$. 

Without loss of generality we may assume the chosen representative 
of a conjugacy class in $\gru{SU}{2}$ is of the form 
$g(\psi)=\exp(i{\psi\over 2}\sigma_3)$ with $0<\psi<2\pi$ (excluding trivial cases), so that its centralizer $Z_\psi$ 
is the circle group  $\gru{U}{1}=
\{\exp(it\sigma_3)\bigm|t\in{\bf R}\}\subset \gru{SU}{2}$.   
Recall the surjective map
\begin{equation}
\label{cover}
{\bf R}^3\ni{\bf x}\mapsto \exp
(i {\bf x}\cdot {\bf \sigma} )\in \gru{SU}{2} \label{**}
\end{equation}
where 
$ {\bf x}\cdot {\bf \sigma} =(x_1\sigma_1+x_2\sigma_2+x_3\sigma_3)\in \ali{su}{2}$ 
and the adjoint representation $g\mapsto \mbox{\rm Ad}(g)$ defined 
by means of the relation 
\[ 
g \exp (i {\bf x}\cdot {\bf \sigma} )g^{-1}
=\exp (i  \mbox{\rm Ad}(g){\bf x}\cdot {\bf \sigma} ), 
\]
which, when the matrices of  $\mbox{\rm Ad}(g)$ are taken with 
respect to the basis $\{i\sigma_\alpha\}_{\alpha=1}^3$, gives the
standard covering $\gru{SU}{2}\to \gru{SO}{3}$.  One knows the map (\ref{cover}) is injective for $|{\bf x}|<\pi$. 

Now expressing the map  
$\gru{SU}{2}/Z_\psi\ni gZ_\psi\mapsto gg(\psi)g^{-1}\in \gru{SU}{2}$
 in terms of the coordinatization given by (\ref{cover}) one arrives at the identification  
of the conjugacy class  $C_\psi$ of  
$g(\psi)$ with the sphere $S_\psi=\{i{\psi\over2}{\bf n}\cdot\sigma
\mid {\bf n}\in {\bf R}^3,\;|{\bf n}|=1\}$ of radius $|{\psi\over2}|$ in ${\bf R}^3\simeq\ali{su}{2}$. 
To see this consider  $G\ni g\mapsto i\mbox{\rm Ad}(g)\sigma_3\in \ali{su}{2}$ --- since the map is constant on $Z_\psi$ and has the unit sphere $S\subset  {\bf R}^3\simeq\ali{su}{2}$ as its image, it gives rise to a bijection of $G/Z_\psi$ with $S$. 
In fact, if we parametrize $\gru{SU}{2}$ by the Euler angles, 
\[
g=g(\varphi,\theta,\psi)=g(\varphi)h(\theta)g(\psi), 
\]
with $g(\varphi),\,g(\psi)$ as above and $h(\theta)=\exp (i {\theta\over2}\sigma_1)$, then we have   
\[
\mbox{\rm Ad}(g(\varphi,\theta,\psi))\sigma_3={\bf n}(\varphi,\theta)\cdot \sigma, 
\]
where
\[
{\bf n}(\varphi,\theta)=
(\sin\theta\sin\varphi,\,\sin\theta\cos\varphi,\,\cos\theta), 
\]
so that $[0,2\pi[\times[0,\pi]\,\ni(\varphi,\theta)\mapsto
\exp(i{\psi\over2}{\bf n}(\varphi,\theta)\cdot \sigma)$ is a
parametrization of  $C_\psi$. 
Since the normalized Haar measure on $\gru{SU}{2}$ is expressed in terms 
of the Euler angles by 
\[
\int_Gf(g)\,dg={1\over16\pi^2}\int_{-2\pi}^{2\pi}\int_0^\pi\int_0^{2\pi}f(g(\varphi,\theta,\psi))\,
\sin\theta\,d\varphi\,d\theta \,d\psi 
\]
it is clear that the invariant integral on $C_\psi$ 
defined by (\ref{InvInt}) is given as 
\[
{1\over4\pi}\int_0^\pi\int_0^{2\pi}f\Bigl(\exp(i{\psi\over2}
{\bf n}(\varphi,\theta)\cdot \sigma)\Bigr)\sin\theta\,d\varphi\,d\theta . 
\]
The class operator for a representation $(T,V)$ will therefore be given by the integral 
\begin{equation}\label{classforSU2}
T(1;g(\psi))= {1\over4\pi}\int_0^\pi\int_0^{2\pi}T\Bigl(\exp(i{\psi\over2}
{\bf n}(\varphi,\theta)\cdot \sigma)\Bigr)\sin\theta\,d\varphi\,d\theta 
\end{equation}
where 
$$
T\Bigl(\exp(i{\psi\over2}
{\bf n}(\varphi,\theta)\cdot \sigma)\Bigr)= T\Bigl(g(\varphi)h(\theta)g(\psi)h(\theta)^{-1}g(\varphi)^{-1}\Bigr).
$$
\end{Example}

\section{The class operator for the group algebra of  a compact group}
In this section we present a rigorous construction of the class operator for the (left) regular representation  (the group algebra) of  compact groups. 
To motivate the subsequent considerations  we  start with 
a brief overview of the construction of the class operator in the case of the finite group, but before doing so,  let us  first introduce some general notations and recall few known facts. 

If $G$ is a compact group we shall always  assume the Haar measure has been normalized so that $\int_G\,dg =1$.  $\widehat{G}$  will denote the set of equivalence classes of irreducible representations of $G$ and for any $\alpha \in \widehat{G}$  $n^\alpha$ will stand for its dimension. 
Given  $\alpha \in \widehat{G}$  we denote $(T^\alpha,\, V^\alpha)$  any of its representatives and 
let  $t^\alpha_{ij}(g)$ be the matrix elements of $T^\alpha(g)$, which we shall assume to satisfy 
$t^\alpha_{ij}(g^{-1})=\overline{t^\alpha_{ji}(g)}$, the bar denoting the complex conjugation. As usual $\chi^\alpha$ will denote the character corresponding to the class $\alpha$. 

Recall also that for any unitary representation $(U,\, W)$ of $G$ and any $\alpha \in \widehat{G}$ there is a uniquely determined invariant subspace $W^\alpha \subset W$, with the property that the restriction of $U$ to $W^\alpha$ is a multiple of $T^\alpha$ --- such a subspace or better the representation associated to it is called the isotypic component of type $\alpha$ of the representation $(U,\, W)$. The decomposition of $(U,\, W)$ into the direct sum of isotypic components, usually called the canonical decomposition, is unique up to an order of summands and is writen in the form 
\[ 
W={\textstyle\bigoplus\limits_{\alpha\in\widehat{G}}}W^\alpha, 
\]
with the corresponding orthogonal projections on $W^\alpha$ given by 
\[ 
 P^\alpha= n^\alpha \int_G\overline\chi^\alpha(g)U(g)\,dg.
\]

\subsection{A special case ---  finite groups}
Assume now  $G$ to be finite; for any subset $X\subset G$ the number of elements of $X$ will be denoted by $| X |$. 
Recall that in the finite case the group algebra ${\cal A}(G)$ is spanned by group elements $g\in G$, assumed to be linearly independent, and the multiplication in ${\cal A}(G)$ is obtained by 
extending the group multiplication from basis elements to the whole algebra by bilinearity. 
The inner product  in ${\cal A}(G)$  is obtained by declaring  the group elements to be mutually orthogonal and their norms set to be equal $| G |^{-1/2}$. Hence if $\phi =\sum_{g\in G}\phi(g)g$,
$\psi =\sum_{g\in G}\psi(g)g$ are arbitrary elements of  ${\cal A}(G)$ then we have
\begin{equation}
(\phi\mid \psi)= \dfrac{1}{| G |}\sum_{g\in G}\phi(g)\overline\psi(g).
\label{inner} 
\end{equation}
It is sometimes convenient not to distinguish between the element  $\phi =\sum_{g\in G}\phi(g)g$ of the group algebra  ${\cal A}(G)$ and the coefficient function $g\mapsto \phi(g)$ describing its coefficients with respect to the basis consisting of the
group elements. Note however that this is in fact an identification of  ${\cal A}(G)$ with 
${\cal K}(G)$ of the previous section, or in yet another (and more common in this case, cf. \cite{WM}) notation with ${\cal F}(G)$, the latter denoting the space of all functions on $G$. 

Let now $C_0\subset G$ be a conjugacy class, $g_0\in C_0$ an arbitrary element and $Z_0$ its centralizer --- recall $|G|=|Z_0||C_0| $. Now let 
\begin{equation}
L_0= \dfrac{1}{|C_0|}\sum_{g\in C_0}g =\dfrac{1}{|G|}\sum_{g\in G}gg_0g^{-1}
\label{finclass} 
\end{equation}
be the class sum corresponding to $C_0 $. If $\chi$ is a function on $G$ constant on conjugacy classes, in particular the character of a representation of $G$, then 
\begin{equation}
( L_0\mid \chi)= \dfrac{1}{|G||C_0|}\sum_{g\in C_0}\overline\chi(g)=
 \dfrac{1}{|G|} \overline\chi(C_0),
\label{innclass}
\end{equation}
where by a slight abuse of notation we have used $\overline\chi(C_0 )$  to denote the common value of  $\chi$ on the elements of the class $C_0$. As $L_0$ clearly belongs to the center of the group algebra  ${\cal A}(G)$  and recalling that irreducible characters form an orthonormal (w. r. to the inner product (\ref{inner})) basis of the centre of ${\cal A}(G)$  we see that
\begin{equation}
L_0= \dfrac{1}{|G|}\sum_{\alpha \in \widehat{G}}\overline\chi^\alpha(C_0)\chi^\alpha.
\label{sum_for_class}  
\end{equation}

Since the left multiplication in $ {\cal A}(G) $ is just the linear extention of the left regular representation $\lambda$, the operator of the left multiplication with $\phi\in  {\cal A}(G)$ will be denoted by $\lambda(\phi)$ and since this action of ${\cal A}(G)$ on itself is faithful,  
we can identify  $\phi$ with   $\lambda(\phi)\in L({\cal A}(G))$. 
In terms of the coefficient function  the left multiplication by $\phi\in  {\cal A}(G)$ corresponds to the convolution  with the function $\widetilde\phi= |G|\phi$, 
i.e.,  if  $\psi =\sum_{g\in G}\psi(g)g \in {\cal A}(G)$,  then  
\begin{equation} 
\lambda(\phi)\psi=\phi\cdot\psi =
\sum_{x\in G} \dfrac{1}{|G|}\sum_{g\in G} \widetilde\phi(g)\psi(g^{-1}x)   x=
\sum_{x\in G} \widetilde\phi \ast\psi(x) x,
\end{equation}
where as usual  the convolution of  $\kappa, \rho\in  {\cal A}(G) $ is denoted by  $\kappa  \ast\rho$  
and defined by the following equality  
\[
\kappa  \ast\rho(x)=\dfrac{1}{|G|}\sum_{g\in G} \kappa(g)\rho(g^{-1}x).
\]
In the effect  
the operator of multiplication with the class sum $L_0$,  which is  the class operator  to be denoted  $\lambda(L_0)=\lambda(L_0;g_0)$,   can be written in the form
 \begin{equation} 
\label{sum_for_class2} 
\lambda(L_0;g_0) = \sum_{\alpha \in \widehat{G}}\dfrac{1}{n^\alpha}\chi^\alpha(C_0)P^\alpha,
\end{equation}
where  $P^\alpha$ is  given by the convolution operator
$\psi\mapsto P^\alpha\psi= n^\alpha \overline\chi^\alpha\ast\psi$.
We shall see presently that the equality (\ref{sum_for_class2})  retains literally its form also for non finite compact groups,  but the above proof, based on the summation over conjugacy class, cannot be adapted to that case.

\subsection{The general  case}
Let now $G$ be an arbitrary compact group. 
In this case  as the group algebra the space $L^2(G) $  of all (equivalence classes of) square integrable functions on $G$  is taken with the multiplication given by the convolution
\[
\phi  \ast\psi(x)=\int_{G} \phi(g)\psi(g^{-1}x)\,dg
\]
and the usual $L^2(G) $  inner product 
\[
(\phi\mid \psi)=\int_{G}  \phi(g)\overline\psi(g)\,dg.
\]
The group $G$ acts on $L^2(G) $ unitarily by means of both left and right regular representations 
(cf. Eq. (\ref{regular})) and again as in the finite case the left or right multiplication (i.e. convolution) can be regarded as the natural extention of the corresponding regular representation of $G$ .   
The isotypic decomposition of the regular representation has the following  form. The spaces  
 $L^2(G) ^\alpha={\rm span}\{\overline{t^\alpha_{ij}}(g)\bigm| 1\le i,j\le n^\alpha\}\subset  L^2(G)$  are minimal subspaces of $L^2(G) $  invariant under left and right regular representation of $G$
and the restriction of $\lambda$ to $L^2(G) ^\alpha$ is an $n^\alpha$-fold multiple of the representation of the class $\alpha$. Moreover, the normalized matrix elements $(n^\alpha)^{1/2}\overline{t^\alpha_{ij}}$ form an orthonormal basis of $L^2(G)^\alpha$. 

The class operator is now, according to (\ref{int2}), given by the integral 
\[
\widetilde{\lambda}(1;g_0)=\int_{G/Z_0 }\lambda(\dot{x})\,d\mu(\dot{x}) =
 \int_{G}\lambda(xg_0x^{-1})\,dx   
\]
 the last equality following in virtue of  compactness of $G$.  We denote it for brevity by $\lambda(C_0)$. The integral has to be evaluated pointwise for $f\in  L^2(G)$, so that 
\begin{equation}\label{classlambda}
\lambda(C_0)f(y)=  \int_{G}\lambda(xg_0x^{-1})f(y)\,dx =  \int_{G}f(xg_0^{-1}x^{-1}y)\,dx.
\end{equation}
We compute the integral on each of the isotypic subspaces separately, so taking  $f=\overline{t^\alpha_{ij}}$ we use the multiplicativity of the matrix elements to get 
\[
\overline{t^\alpha_{ij}}(xg_0^{-1}x^{-1}y)= \sum_{k=1}^{n^\alpha}\overline{t^\alpha_{ik}}(xg_0^{-1}x^{-1})\overline{t^\alpha_{kj}}(y).
\]
Inserting this expression into (\ref{classlambda}) we obtain
\[
\lambda(C_0) \overline{t^\alpha_{ij}}(y)=\sum_{k=1}^{n^\alpha} \overline{t^\alpha_{kj}}(y)
 \int_{G}\overline{t^\alpha_{ik}}(xg_0^{-1}x^{-1})\,dx.
\]
Expanding further and using the orthogonaliy of matrix elements we get after easy manipulations 
\[
 \int_{G}\overline{t^\alpha_{ik}}(xg_0^{-1}x^{-1})\,dx=
\dfrac{1}{n^\alpha}\overline{\chi^\alpha}(g_0^{-1})\delta_{ik} = \dfrac{1}{n^\alpha}\chi^\alpha(g_0)\delta_{ik}.
\]
Thus $\lambda(C_0)$ acts on $L^2(G) ^\alpha$ as multiplication by $(n^\alpha)^{-1}\chi^\alpha(g_0)$, 
so we have obtained the result claimed above.
\begin{Prop}
The class operator for the group algebra of a compact group $G$ based on a conjugacy class $C_0$ is given by the formula 
\[
\lambda(C_0)=\sum_{\alpha \in \widehat{G}}\dfrac{1}{n^\alpha}\chi^\alpha(C_0)P^\alpha,
\]
with $\chi^\alpha(C_0)$ denoting the (common) value of $\chi^\alpha$ on (the elements of) the class 
$C_0$. 
\end{Prop}
We note in particular, that for the case of the group  $\gru{SU}{2}$ and the class operator given by the formula (\ref{classforSU2})  with $T$ replaced by $\lambda$,  this gives the following expression (cf. \cite{FR,B,OS,MZ})
\[
\lambda(C_\psi)= \sum_{j=0}^{\infty}\dfrac{\sin(2j+1){\psi\over2}}{(2j+1)\sin{\psi\over2}}P^j,
\]
with $P^j$ denoting projections on subspaces corresponding to the spin $j$. 

\section{Connection with tensor operators}
In this last section we consider  the connections of the above construction of weighted class operators with the concept of tensor operators which is known as being of  paramount importance for physics. The general discussion, which was given in a recently published paper of A. Or{\l}owski and the author \cite{SO}, will be briefly recalled first and then suplemented by the  computation of weighted class operators for the group algebra along the lines of section 3.2. As above, the discussion will  be confined to the compact case (for some earlier related work, see \cite{dV-Z,K-D,B-G}).  

\subsection{Definitions and general results on tensor operators}
In a conventional formulation a tensor operator is understood to be 
a linearly independent set of operators $\{T_i\}$ on the vector space 
of a certain representation $(U,\,H)$ of $G$, such that 
\begin{equation} 
\label{TOI} 
U(g) T_i U(g)^{-1}= \sum_{j}D_{ji}(g)T_j, 
\end{equation} 
where  $D_{ij}(g)$ are matrix coefficients of a representation of the group $G$. 

However,  recall the following   
elegant (although purely algebraic) definition given by  L.~Michel in \cite{LM}. 

\begin{Definition} Let $(S,\, V)$ and  $(U,\, H)$ be two  
representations  of a group $G$. Then a tensor operator of type $(S,\,V)$ over 
the space $H$ is a (nonzero) linear mapping  $T:V\to L(H)$,  intertwining $S$ 
with the conjugation representation defined  by $U$  on $L(H)$. 
\end{Definition} 
Since the intertwining condition means that the following equality is valid 
\begin{equation} 
\label{TO} 
T (S(g)v)= U(g) T(v) U(g)^{-1},\qquad 
\forall\, g\in G,\ \forall\, v\in V, 
\end{equation} 
it is clear that the results of Section 2  are in fact statements 
about tensor operators. To describe then properly we shall need some more terminology.

A tensor operator $T:V\to L(H)$ is called irreducible if the 
representation $(S,\,V)$ is irreducible. If $\alpha $ is the class of 
representation $(S,\,V)$, then we say that $T$ is of the type $\alpha $. 
Two irreducible tensor operators of the same type $T_1:V_1\to L(H)$ and 
$T_2:V_2\to L(H)$ are called independent if their images are different 
subspaces of $L(H)$ --- note that the images being irreducible they can 
either be identical or have only zero in common. 

Before going further, let us indicate that the equivalence of these two notions of tensor operators can be obtained by  choosing  a basis, say  $\{v_i\}$, of the space $V$ and setting  
$T_i=T(v_i)$. Then (\ref{TOI}) follows trivially from (\ref{TO}) by taking $D_{ij}(g)$
 to be the matrix coefficients of $S(g)$ defined by the usual recipe
$S(g)v_i=\sum_{j=1}^{\dim{V}}D_{ji}(g)v_j$. 

Now  we can state 
\begin{Theo}{\rm(cf. \cite{SO})} \label{ }
The weighted class operator map 
{\rm(\ref{def2})} is a tensor operator of the type of the induced representation 
$(\lambda,\,L^2(G/Z_0))$. The restriction of this map to any 
invariant and irreducible subspace $V^\alpha  \subset  L^2(G/Z_0)^\alpha $, if nonzero, is an 
irreducible tensor operator of type $\alpha $. 
\end{Theo}

In order to actually construct irreducible tensor operators as weighted class operators
one needs to know explicitely  the canonical decomposition of the induced representation, 
or better still,  irreducible subspaces of $L^2(G/Z_0)$.  The first part of the issue is 
solved by the classical Frobenius reciprocity theory, cf.  e.g. Section 4.3 of \cite{AW}
or Chapter 8 of  \cite{AR}.  

Let  $\widehat{G}_{0}\subset \widehat{G}$ be the 
subset consisting of (classes of) representations admitting nonzero 
${Z_0}$-fixed vectors. Given $\alpha \in \widehat{G}_{0}$ choose $(T^\alpha,\,V^\alpha)\in  \alpha$ and let $V^\alpha_0\subset V^\alpha $ be the subspace of ${Z_0}$-fixed vectors. Set 
$m^\alpha=\dim V^\alpha_0$ and note this number depends on the class $\alpha$, but not on the choice of  $(T^\alpha,\,V^\alpha)\in  \alpha$. On the other hand consider the induced representation  $(\lambda,\, L^2(G/Z_0))$ and let ${\bf i}(T^\alpha\colon \lambda)$ be the multiplicity of the representation  $T^\alpha$ in the induced representation. Then the Frobenius reciprocity states these two numbers are equal, $m^\alpha={\bf i}(T^\alpha\colon \lambda)$. 

Consequently the subspace  $L^2(G/Z_0)^\alpha $  is nonzero if and only if $\alpha \in \widehat{G}_{0}$ and in this case it has dimension $m^\alpha n^\alpha $  and contains  $m^\alpha $ copies of  representations  of the class $\alpha$. 
Thus the only irreducible tensor operators  which can  be constructed as weighted class operators based on the  conjugacy class $C_0=G/Z_0$ are those  whose types  contain the trivial (identity) representation of $Z_0$. 

\subsection{Tensor operators over the group algebra}

We shall now give a description of the construction of  tensor operators over the group algebra extending the method of Section 3.2.  The aim is to compute the integrals 
\begin{equation} \label{tensorlambda}
\widetilde{\lambda}(\widetilde{f};\,g_0)=
\int_{G/Z_0 }\widetilde{f}(\dot{x})\lambda(\dot{x})\,d\mu(\dot{x}) =
 \int_{G}f(x)\lambda(xg_0x^{-1})\,dx,  
\end{equation}
where  $f$ and $\widetilde {f}$ are related by (\ref{transfer}). Since we are interested in irreducible tensor operators we can assume that $\widetilde{f} $'s belong to an irreducible subspace of  $L^2(G/Z_0)^\alpha $. Such functions can be constructed in the following way. 
Choose  an orthonormal basis $\{e_i\}$ in  $V^\alpha$ in such a way 
that its first $m^\alpha $ vectors form a basis in  $V^\alpha_0$ and 
the remaining ones span the complementary $Z_0$-invariant 
subspace and denote by $t^\alpha_{ij}(g)=\langle T^\alpha(g){e_j}\mid {e_i} \rangle$ 
the matrix elements of  $(T^\alpha,\,V^\alpha)$ with respect to this chosen basis 
(we are using here the inner product which is linear in the first variable).   Note that the matrix elements  $t^\alpha_{ij}(g)$ for $1\le j\le m^\alpha$ are right $Z_0$-invariant functions, hence  can be regarded  as functions from  $L^2(G/Z_0)$ and thus  $t^\alpha_{ij}$ can be identified with their right $Z_0$ averages $\widetilde{t^\alpha_{ij}}$. 
In virtue of the  Frobenius reciprocity  their complex conjugates  $ \overline{t^\alpha_{ij}}(g)$ span the space $L^2(G/Z_0)^\alpha $.    

It follows that for a given representation $(T,\,V)$ of $G$  by setting $T^{(\alpha,j)}_i=\widetilde{T}(\overline{t^\alpha _{ij}};\,g_0)$ we obtain  in general 
 $m^\alpha $ sets  $\{T^{(\alpha ,j)}_i\mid1\le i\le n^\alpha  \}$ 
of weighted class operators  transforming according to the formula 
$V(g) T^{(\alpha ,j)}_i V(g)^{-1}= \sum_{k}t^\alpha _{ki}(g)T^{(\alpha ,j)}_k$. 
However, some of these sets may degenerate (i.e. consists only of  the  zero operator),  
but the author is unaware of  any general criterion for nonvanishing of an intertwining operator on a 
given irreducible subspace. We shall not pursue this question here.

Now, let us substitute a matrix element  $\overline{t^\alpha_{kl}}$ for $f$ in the integral (\ref{tensorlambda}) and evaluate the restriction of $\widetilde{\lambda}(\overline{t^\alpha_{kl}};\,g_0)$ 
to the subspace $L^2(G)^\sigma$. Proceeding exactly as above in the Section 3.2 we get  
\begin{eqnarray} \label{W-E}
\widetilde{\lambda}(\overline{t^\alpha_{kl}} ;\,g_0)\overline{t^\sigma_{ij}} (y)&=&
\int_{G}\overline{t^\alpha_{kl}}(x)\overline{t^\sigma_{ij}} (xg_0^{-1}x^{-1}y)\,dx \nonumber \\
&=& \sum_{s} \overline{t^\sigma_{sj}}(y)
\int_{G}\overline{t^\alpha_{kl}}(x) \overline{t^\sigma_{is}}(xg_0^{-1}x^{-1})\,dx\nonumber \\
&=& \sum_{spr} \overline{t^\sigma_{sj}}(y) t^\sigma_{pr}(g_0) 
\int_{G}\overline{t^\alpha_{kl}}(x)
\overline{t^\sigma_{ir}}(x)t^\sigma_{sp}(x)\,dx.
\end{eqnarray}

Now, by the general wisdom of the canonical decomposition, one knows that the products of matrix elements of irreducible representations can be written as linear combinations of matrix elements of irreducible representations. One way of doing so is to use the so called coupling coefficients (cf. \cite{K-D,SO}).  We digress  briefly to introduce the needed notions, in the form familiar from  the treatment of the classical  Clebsch--Gor\-dan coefficients.   

For $(T^\sigma,\,V^\sigma)\in \sigma$ consider the conjugation representation in $L(V^\sigma)$ and  write its canonical decomposition in the form
\begin{equation}
L(V^\sigma)={\textstyle\bigoplus\limits_{\gamma\in \Gamma(V^\sigma)}}E^\gamma=
{\textstyle\bigoplus\limits_{\gamma\in \Gamma(V^\sigma)}}m(\sigma;\gamma)V^\gamma,
\label{decomp}
\end{equation}
where $\Gamma(V^\sigma)\subset \widehat{G}$ denotes the set of classes
of irreducible unitary representations of $G$ which occur in that decomposition,
$E^\gamma$ is the isotypic component of the class $\gamma$, and the right hand side of the equality is obtained by further decomposing  $E^\gamma$ into irreducibles. The multiplicity  $m(\sigma;\gamma)$ of the class $\gamma$ in $L(V^\sigma)$ is sometimes called $3j$ symbol and denoted $\{\sigma\overline{\sigma}\gamma\}$, cf. e. g. \cite{dV-Z}. 
Following the general usage we are 
using $\overline{\sigma}$ to denote the complex conjugate representation, i.e. the one with the complex conjugate matrix elements.   An important and much simpler situation occurs, when the canonical decomposition is multiplicity free, i.e. when $m(\sigma;\gamma)=1 $ for each $\gamma\in \Gamma(V^\sigma) $, what is the case in particular for the so called simply reducible groups of Wigner.

 The coupling coefficients are related to a choice of a basis realizing the canonical decomposition.  The  vectors of such (orthonormal)  basis are denoted $e^\gamma_{mn}$, where $\gamma$ describes the classes of irreducible representations occurring in the decomposition,
$1\le m\le m(\sigma;\gamma)$ distinguishes between different copies of the same
representation of the class $\gamma$ and $1\le n\le n^\gamma$ indexes vectors
of a given base within a fixed copy of the representation space.  In
particular for fixed $\gamma$ and $m$ the set $\{e^\gamma_{mn}\mid1\le n\le
n^\gamma\}$ is an orthonormal basis for an invariant subspace, on which the
conjugation representation acts by a representation belonging to the class
$\gamma$.  Now fix an orthonormal basis $\{v_i\}$ for $V^\sigma$ and let 
$E_{ij}\in L(V^\sigma)$ be a ``matrix unit" corresponding to that basis, 
i.e., the linear map given by $v\mapsto E_{ij}(v)=\langle{v}\mid{v_j}\rangle{v_i}$. 
The coupling coefficients $c(\sigma i;\overline{\sigma}j\bigm|\gamma{mn})$ 
give the transition between the two bases  (the summation extending over the whole range of indices involved)
\begin{equation} \label{coupling}
E_{ij}=\sum_{\gamma{mn}}c(\sigma i;\overline{\sigma}j\bigm|\gamma{mn})e^\gamma_{mn},
\end{equation} 
and it is a simple exercise to show validity of the equation 
\begin{equation} 
\label{matrixproduct}
\overline{t^\sigma_{ir}}(x)t^\sigma_{sp}(x)=
\sum_{\gamma{ m n q}}\overline{c(\sigma s;\overline{\sigma}i\bigm |\gamma{mq})}t^\gamma_{qn}(x)
c({\sigma} p;\overline{\sigma}r\bigm|\gamma{mn}).
\end{equation}
It now follows that 
\begin{eqnarray*} 
\int_{G}\overline{t^\alpha_{kl}}(x)
\overline{t^\sigma_{ir}}(x)t^\sigma_{sp}(x)\,dx& =&
\sum_{\gamma{ m n q}}\overline{c(\sigma s;\overline{\sigma}i\bigm |\gamma{mq})}
c({\sigma} p;\overline{\sigma}r\bigm|\gamma{mn})
\int_{G}\overline{t^\alpha_{kl}}(x) t^\gamma_{qn}(x)\,dx\\
&=&\dfrac{1}{n^\alpha}\sum_{m}\overline{c(\sigma s;\overline{\sigma}i\bigm |\alpha{mk})}
 c({\sigma} p;\overline{\sigma}r\bigm|\alpha{ml}).
 \end{eqnarray*}
Inserting this into the integral (\ref{W-E}) we get 
\begin{eqnarray}
\widetilde{\lambda}(\overline{t^\alpha_{kl}} ;\,g_0)\overline{t^\sigma_{ij}} (y) & = & 
\dfrac{1}{n^\alpha} \sum_{mspr} \overline{t^\sigma_{sj}}(y) t^\sigma_{pr}(g_0) 
\overline{c(\sigma s;\overline{\sigma}i\bigm |\alpha{mk})}
 c({\sigma} p;\overline{\sigma}r\bigm|\alpha{ml}) \nonumber \\
 &= & \dfrac{1}{n^\alpha} \sum_{sm}\overline{t^\sigma_{sj}}(y)
\overline{c(\sigma s;\overline{\sigma}i\bigm |\alpha{mk})}
\biggl( \sum_{pr} c({\sigma} p;\overline{\sigma}r\bigm|\alpha{ml}) t^\sigma_{pr}(g_0)\biggr).  
\end{eqnarray}
From that one immediately computes the matrix coefficients of the operator  $\widetilde{\lambda}(\overline{t^\alpha_{kl}} ;\,g_0)$  obtaining this way  the following result of the Wigner--Eckart type.
 
\begin{Theo} Let for  $\sigma, \gamma\in \widehat{G}$  the coupling coefficients  
$c(\sigma i;\overline{\sigma}j\bigm|\gamma{mn})$ be defined by the equation {\rm(\ref{coupling})}.
For any matrix coefficient $\{t^\alpha_{kl}\}$ of the irreducible representation of $G$ of the class  $\alpha$  let  $\widetilde{\lambda}(\overline{t^\alpha_{kl}} ;\,g_0)$ be the corresponding  weighted class operator  for the group algebra $G$ defined by the equation {\rm(\ref{tensorlambda})}. 
Then its matrix coefficients are given by the formulae of the form 
\begin{eqnarray} \label{W-E2} 
n^\sigma\langle\widetilde{\lambda}(\overline{t^\alpha_{kl}} ;\,g_0) \overline{t^\sigma_{ij}}
\mid \overline{t^\gamma_{uv}} \rangle & =& \dfrac{1}{n^\alpha}\delta_{\gamma\sigma}\delta_{jv} 
\sum_{m=1}^{m(\sigma;\alpha)}\overline{c(\sigma u;\overline{\sigma}i\bigm |\alpha{mk})}
\nonumber \\
&& \times 
\biggl( \sum_{pr} c({\sigma} p;\overline{\sigma}r\bigm|\alpha{ml}) t^\sigma_{pr}(g_0)\biggr)
\end{eqnarray} 
\end{Theo} 
Let us point out that in addition to the usual content  of the Wigner--Eckart theorem we have obtained in  the formula  (\ref{W-E2})  an explicit and complete description of the reduced matrix elements  which are given by  
$$
 \dfrac{1}{n^\alpha} \sum_{pr} c({\sigma} p;\overline{\sigma}r\bigm|\alpha{ml}) t^\sigma_{pr}(g_0)
$$
depending on the chosen element $g_0$ determining the conjugacy class. 

\section*{Acknowledgments}
{\footnotesize
A preliminary version of the results contained in the paper  was presented during the  
II~International Workshop  "Lie Theory and Its Applications in Physics"  at the Arnold-Sommerfeld Institute of the Technical University Clausthal.  Thanks are  due to the organizers of the Workshop  for the invitation to present these results and for creating so pleasant and inspiring  atmosphere of the meeting. 

The work on this paper was  partially supported by State Committee for Scientific  Research (SCSR) grant No 2 P03B 044 13, and also by the research project SFB 288 --- {\em Differentialgeometrie und Quantenphysik} at the Institut f\"ur Reine  Mathematik of the Humboldt University in Berlin,  where the actual writing of the paper took place. The author expresses his thanks to Prof. Thomas Friedrich for the invitation and hospitality during his stay in Berlin.
}

\end{document}